\documentclass[a4paper, 11pt, twoside]{article}
\usepackage{amsmath, amssymb}
\usepackage{doc, exscale, fontenc, latexsym, syntonly}
\usepackage{xypic}

\author{ A. Mutlu, T. Porter}

\title{\huge{Iterated Peiffer pairings in the Moore complex \\
 of a simplicial group}}

\newtheorem{defn}{Definition}[section]
\newtheorem{prop}[defn]{Proposition}
\newtheorem{thm}[defn]{Theorem}
\newtheorem{lem}[defn]{Lemma}
\newtheorem{cor}[defn]{Corollary}

\newenvironment{pf}{{\bf Proof:}}{$\Box$\mbox{}}

\begin{document}

\maketitle
\begin{abstract}
We introduce a pairing structure within the Moore complex ${\bf NG}$ of a simplicial group ${\bf G}$ and use it 
to investigate generators for $NG_n\cap D_n$ where $D_n$ is the subgroup generated by degenerate elements. This
is applied to the study of algebraic models for homotopy types.\\
{\noindent\bf A. M. S. Classification:}  18D35 18G30 18G50 18G55.
\end{abstract}
\section*{Introduction}
Simplicial groups and simplicial groupoids are valuable algebraic models for homotopy types.
Much has been studied about the way the group structure interacts with the simplicial structure 
to yield homotopy information. 

Recently the work of Wu, \cite{wut}, ~\cite{wu}, has shown that there is still progress that 
can be made in calculation of homotopy invariants such as homotopy groups from simplicial groups.
Wu used techniques of combinatorial group theory, iterated commutators and properties related 
to the semidirect product decompositions of the individual $G_n$ to give some insight into, for
instance, $\pi_{n+1}(\Sigma K(\pi,1))$, the homotopy groups of the suspension of an Eilenberg-MacLane 
space.

Earlier Brown and Loday, \cite{bl1}, had used techniques derived from their generalised van Kampen 
theorem and Loday's theory of cat\textsuperscript{$n$}-groups to give a complete description of the 
$3$-type of $\Sigma K(\pi,1).$ This raises the possibility of linking the results of Wu with crossed 
algebraic techniques and to combine the two techniques in order to give descriptions of, for instance,
the $k$-type of $\Sigma K(\pi,1)$ for $k=4$ and $5.$ This is still out of our 
reach with the techniques of this paper, but other results 
suggest the way to develop tools for this sort of task. 

Carrasco \cite{carrasco}, and with Cegarra in
\cite{c:c}, gave a complete description of the extra structure of the Moore complex, ${\bf NG},$ of a 
simplicial group ${\bf G}$ needed to reconstruct ${\bf G}$ from ${\bf NG}$, a sort of ultimate 
generalisation of the classical Dold-Kan theorem that links simplicial abelian groups with chain 
complexes. The controlled vanishing of this extra structure given necessary and sufficient conditions 
for the Moore complex to be a crossed complex or crossed chain complex. Further links between simplicial 
groups, their Moore complexes and crossed algebraic models for homotopy types have been given by Baues 
\cite{baues1},~\cite{baues2} and \cite{baues3} and also by the second author \cite{porter}.

In this article we will develop a variant of the Carrasco - Cegarra pairing operators, that we 
will call {\it Peiffer pairings}, and will show that these pairings give products of commutators, and thus, 
by repeated application, iterated commutators that generate the Moore complex terms in those dimensions 
where additional non-degenerate generators are not present and in general, they generate $NG_n\cap D_n$ 
where $D_n$ is the subgroup of $G_n$ generated by the degenerate
elements. So far it has not been possible to find a general form for the
relations between these generators.  This would seem to be an extremely
hard problem in general.  Some results in low dimensions and for free
simplicial groups have been obtained, but as they are incomplete they
will not be included here.
Some sample calculations of these generating elements will 
be given as will some fairly elementary examples of the use of this result. \\
{\noindent\bf Acknowledgement}\\
A. Mutlu wishes to thank the University of Celal Bayar, Manisa, Republic of Turkey, for the award of a 
research scholarship during the tenure of which this work was undertaken.
\section{Simplicial groups, Moore complexes and Peiffer pairings}
We refer the reader to Curtis's survey article \cite{curtis} or May's book, 
\cite{may}, for most of the basic properties of 
simplicial sets, simplicial groups, etc. that we will be needing.
\subsection{The Moore complex}
If ${\bf G}$ is a simplicial group, the Moore complex $({\bf NG},\partial)$ of ${\bf G}$ is the 
(non-abelian) chain complex defined by 
$$
NG_n = \bigcap\limits_{i=0}^{n-1}\text{Ker}d_i
$$
with $\partial_n : NG_n\longrightarrow NG_{n-1}$ induced from $d_n^n$ by restriction. It is well known 
that {\em  n\textsuperscript{th} homotopy group} $\pi_n({\bf G})$ of ${\bf G}$ is the 
{\em  n\textsuperscript{th} homology } of the Moore complex of ${\bf G}$ 
$$
\begin{array}{rcl}
\pi _n({\bf G}) & \cong & H_n( 
{\bf NG},\partial ) \\  & = & \bigcap\limits_{i=0}^n\text{Ker}%
d_i^n/d_{n+1}^{n+1}(\bigcap\limits_{i=0}^n\text{Ker}d_i^{n+1}). 
\end{array}
$$

{\bf Remark and Warning}

There is a possibility of confusion as to the exact definition of $\bf
NG$ as two conventions are currently used, one as above takes the
intersection of the $\text{Ker}d_i$ for $i < n$, the other the
intersection of the $\text{Ker}d_i$ for $0 < i \leq n$. (Curtis
\cite{curtis} uses this latter convention, whilst May, \cite{may}, uses
the former.)  The two theories run parallel and are essentially `dual'
to each other, however there is a necessity for checking, which
convention is being used in any source as the actual form of any formula
usually depends on the convention being used.

\subsection{The poset of surjective maps}
We recall the following notation and terminology referring the reader 
to the work of Conduch\'e, \cite{conduche}, Carrasco and Cegarra \cite{c:c} for more motivation and 
some related results. 

For the ordered set $[n] = \{ 0 < 1 < \cdots < n \}$, let $ {\alpha}_{i}^{n}
: [n+1] \to [n]$ be the increasing surjective  map given by
\[ \alpha_{i}^{n}(j) = \left\{ \begin{array}{ll} 
                         j    &  \mbox{ if $ j \leq i $ } \\
                         j-1  &  \mbox{ if $ j > i $ }
                        \end{array}
                        \right. \]
Let $S(n , n-l)$ be the set of all monotone increasing surjective maps from 
$[n]$ to $[n-l].$ 
This can be generated from the various $ \alpha_{i}^{n}$
by composition. The composition of these generating maps satisfies the
rule $ \alpha_{j} \alpha_{i} = \alpha_{i-1} \alpha_{j}$ with $j < i$. This 
implies that every element $ \alpha \in S(n , n-l)$ has a unique expression 
as $ \alpha = \alpha_{i_{1}} \alpha_{i_{2}} \ldots  
\alpha_{i_{l}}$  with $0 \leq i_{1} < i_{2} < \cdots < i_{l} \leq n$, where 
the indices $i_{k}$ are the elements of $[n]$ at which $\{ i_{1}, ..., i_{l}, 
\} = \{ i : \alpha(i) = \alpha(i+1) \}$. We thus can identify $S(n , n-l)$
with the set $\{(i_{l}, ..., i_{1}) : 0 \leq i_{1} < i_{2} < \cdots < i_{l} 
\leq n-1 \}$. In particular the single element of $S(n ,n)$, defined by the 
identity map on $[n]$, corresponds to the empty $0$-tuple $(~)$ denoted by 
$\emptyset_{n}$. Similarly the only element of $S(n, 0)$ is $(n-1 , n-2 ,  
\ldots , 0)$. For all $n \geq 0$, let 
\[ S(n) = \bigcup_{0 \leq l \leq n} S(n , n-l) .\]
We say that $ \alpha = (i_{l}, ..., i_{1}) > \beta = (j_{m}, ..., 
j_{1})$  in $S(n)$
\[ \mbox{ if $i_{1} = j_{1} , \cdots , i_{k} = j_{k}$ but $i_{k+1} < j_{k+1}$
$( k \geq 0 )$} \]   or
\[ \mbox {if \qquad $i_{1} = j_{1} , \cdots , i_{m} = j_{m}$  and $l > m$}. \]
This makes $S(n)$ an ordered set. For instance, the orders of $S(2)$ and  
$S(3)$ and $S(4)$ are respectively:

$S(2) = \{ \emptyset_{2} < (1) < (0) <(1 , 0) \}$,

$S(3) = \{ \emptyset_{3} < (2) < (1) < (2 , 1) < (0) < (2 , 0) < (1 , 0) 
< (2 , 1 , 0) \}$, 

$S(4) = \{ \emptyset_{4} < (3) < (2) < (3 , 2) < (1) < (3 , 1) < (2 , 1) 
< (3 , 2 , 1) < (0) < (3 , 0) < (2 , 0) < (3 , 2 , 0) < (1 , 0) 
< (3 , 1 , 0) < (2 , 1 , 0) < (3 , 2 , 1 , 0) \}$.  

If $ \alpha, \beta \in S(n),$ we define $ \alpha \cap \beta $ to be the set of 
indices which belong to both $\alpha$ and $\beta.$ 

If $\alpha = (i_{l}, ..., i_{1})$, then we say $\alpha$ has length $l$ and
will write $\# \alpha = l$.
\subsection{The semidirect decomposition of a simplicial group}
The fundamental idea behind this can be found in Conduch\'e \cite{conduche} .
A detailed investigation of this for the case of  simplicial groups is given 
in Carrasco and Cegarra \cite{c:c}.

\begin{lem}\label{1} Let ${\bf G}$ be a simplicial group. Then $G_{n}$ can be decomposed as 
a semidirect product:
\[ G_{n} \cong ~\text{Ker}d_0^{n} \rtimes s_{0}^{n-1}(G_{n-1}) \]
\end{lem}
\begin{pf} The isomorphism can be defined as follows:
\[ \theta : G_{n} \to \text{Ker}d_{0}^{n}\rtimes s_{0}^{n-1}(G_{n-1}) \] 
\[ g \mapsto (g s_{0}d_{0}g^{-1}, s_{0}d_{0}g)  .\]
\hfill\end{pf}\\ 
Since we have the isomorphism $G_{n} \cong ~\text{Ker}d_{0} \rtimes s_{0}G_{n-1}$, 
we can repeat this process as often as necessary to get each of the $G_{n}$ 
as a multiple semidirect product of degeneracies of terms in the Moore 
complex. In fact, let ${\bf K}$ be the simplicial group defined
by 
$$
\begin{array}{cccc}
K_n=\text{Ker}d_{0}^{n+1}, & d_i^n=d_{i+1}^{n+1}\mid _{\text{Ker}
d_{0}^{n+1}} & \text{and} & s_i^n=s_{i+1}^{n+1}\mid _{\text{Ker}
d_{0}^{n+1}}. 
\end{array}
$$
 Applying Lemma ~\ref{1} above, to $G_{n-1}$ and to $K_{n-1}$, gives 
$$
\begin{array}{rcl}
G_n & \cong & \text{Ker}d_0\rtimes s_{0}G_{n-1} \\  
& = & \text{Ker}d_0 \rtimes s_{0}(\text{Ker}d_{0}\rtimes s_{0}G_{n-2}) \\  
& = & K_{n-1} 
\rtimes \ (s_{0}\text{Ker}d_{0}\rtimes  s_{0}s_{0}G_{n-2}). 
\end{array}
$$
Since ${\bf K}$ is a simplicial group, we have the following 
$$
\begin{array}{rcl}
\text{Ker}d_0=K_{n-1} & \cong & \text{Ker}d_{0}^K\rtimes s_{0}^KK_{n-2} \\
& = & (\text{Ker}d_{1}\cap \text{Ker}d_0)\rtimes s_{1}\text{Ker}d_{0} 
\end{array}
$$
and this enables us to write 
$$
G_n=((\text{Ker}d_{1}^n\cap \text{Ker}d_0^n)\rtimes 
s_{1}(\text{Ker}d_{0}^{n-1}))\rtimes
(s_{0}(\text{Ker}d_{0}^{n-1})\rtimes s_{0}s_{0}(G_{n-2})). 
$$
We can thus decompose $G_{n}$ as follows:
\begin{prop}\label{2} (cf. \cite{conduche}, p.158) If ${\bf G}$ is a simplicial group, then for any $n \geq 0$ 
$$
\begin{array}{lll}
G_n & \cong  & (\ldots (NG_n
\rtimes s_{n-1}NG_{n-1})\rtimes \ldots \rtimes s_{n-2}\ldots s_1NG_1)\rtimes
\\& &\quad(\ldots (s_{0}NG_{n-1}\rtimes s_{1}s_{0}NG_{n-2})%
\rtimes \ldots \rtimes s_{n-1}s_{n-2}\dots s_0NG_0).~\square
\end{array}
$$
\end{prop}
The bracketing and the order of terms in this multiple semidirect product
are generated by the sequence:
$$
\begin{array}{lll}
G_1 & \cong & NG_1 
\rtimes s_0NG_0 \\ G_2 & \cong & (NG_2 
\rtimes s_1NG_1)\rtimes (s_0NG_1\rtimes  s_1s_0NG_0) \\ G_3 & \cong & ((NG_3 
\rtimes s_2NG_2)\rtimes (s_1NG_2\rtimes s_2s_1NG_1))\rtimes \\& &\qquad
((s_0NG_2\rtimes s_2s_0NG_1)\rtimes (s_1s_0NG_1\rtimes %
s_2s_1s_0NG_0)). 
\end{array}
$$
and
$$
\begin{array}{lll}
G_4 & \cong & (((NG_4 
\rtimes s_3NG_3)\rtimes (s_2NG_3\rtimes s_3s_2NG_2))\rtimes \\  &  & \qquad
\ ((s_1NG_3 
\rtimes s_3s_1NG_2) \rtimes(s_2s_1NG_2\rtimes s_3s_2s_1NG_1)))\rtimes \\  &  
& \qquad \qquad s_0( \text{decomposition of }G_3). 
\end{array}
$$
Note that the term corresponding to  $\alpha
=(i_l,\ldots ,i_1)\in S(n)$ \ \ is \ \ 
$s_\alpha (NG_{n-\#\alpha
})=s_{i_l...i_1}(NG_{n-\#\alpha })=s_{i_l}...s_{i_1}(NG_{n-\#\alpha }),$
where $\#\alpha =l.$ 
Hence any element $x\in G_n$ can be written in the form%
$$
x=y \prod\limits_{\alpha \in S(n)}s_\alpha (x_\alpha )\text{ \qquad with }%
y\in NG_n\text{ and }x_\alpha \in NG_{n-\#\alpha }. 
$$
\section{Peiffer pairings generate}
In the following we will define a normal subgroup $N_n$ of $G_n.$ First of all we 
adapt ideas from Carrasco \cite{carrasco} to get the construction of a useful
family of natural pairings.
We define a set $P(n)$ consisting of pairs of elements $(\alpha , \beta)$ 
from $S(n)$ with $\alpha \cap \beta =\emptyset $ and $\beta < \alpha $, with respect to 
lexicographic ordering in $S(n)$ where $\alpha =(i_l,\ldots ,i_1),\beta=(j_m,...,j_1)\in S(n).$ 
The pairings that we will need, 
$$
\{F_{\alpha ,\beta }:NG_{n-\#\alpha}\times NG_{n-\#\beta}\longrightarrow
NG_n:(\alpha ,\beta )\in P(n),\ \ n\geq 0\} 
$$
are given as composites by the diagram
$$
\diagram
 NG_{n-\#\alpha} \times NG_{n-\#\beta}   \rto^{\hspace{.9cm} F_{\alpha, \beta}} 
 \dto_{s_{\alpha} \times s_{\beta}} & NG_n \\ 
G_n \times G_n \rto^{\mu}  & G_n  \uto_p                       
\enddiagram
$$
where 
$$
\begin{array}{c}
s_\alpha =s_{i_l}\ldots s_{i_1}:NG_{n-\#\alpha}\longrightarrow G_n, \ \
s_\beta =s_{j_m}\ldots s_{j_1}:NG_{n-\#\beta}\longrightarrow G_n, 
\end{array}
$$
$p:G_n\rightarrow NG_n$ is defined by the composite projections $p(x)=p_{n-1}\ldots
p_0(x),$ where 
$$
p_j(z) = zs_jd_j(z)^{-1}~\qquad \text{with\quad }j=0,1,\ldots, n-1, 
$$
$\mu :G_n\times G_n\rightarrow G_n$ is given by the commutator map and
$\#\alpha$  is  the number of the elements in the set of $\alpha$, similarly for
$\#\beta.$  Thus
$$
\begin{array}{rcl}
F_{\alpha ,\beta }(x_\alpha , y_\beta ) & = & p\mu (s_\alpha \times
s_\beta )(x_\alpha , y_\beta ) \\  
& = & p[s_\alpha x_\alpha   , s_\beta y_\beta ] .  
\end{array}
$$

\textbf{Definition} Let  
 $N_n$ or more exactly $N_n^G$ be the normal subgroup of $G_n$  generated by elements of 
the form 
$$
F_{\alpha ,\beta }(x_\alpha , y_\beta ) 
$$
where $x_\alpha \in NG_{n-\#\alpha}$ ~and ~ $y_\beta \in NG_{n-\#\beta}$.

\medskip

This normal subgroup  $N_n^G$ depends functorially on $G$, but we will usually 
abbreviate  $N_n^G$ to  $N_n$, when no change of group is involved.

\medskip

We illustrate this subgroup for $n=2$ and $n=3$ to show what it looks 
like.\\
\textbf{Example (a) :}
For $n=2,$ suppose $\alpha =(1)$, $\beta =(0)$ and \\
$x_1,y_1\in NG_1={\rm Ker}d_0$. It follows that 
$$
\begin{array}{rcl}
F_{(0)(1)}(x_1 , y_1) & = & p_1p_0[s_0x_1 , s_1y_1] \\  
& = & p_{1}[s_0x_1, s_1y_1] \\  
& = & [s_0x_1, s_1y_1]{~} [s_1y_1, s_1x_1] 
\end{array}
$$
is a generating element of the normal  subgroup $N_2.$

For $n=3,$ the possible pairings are the following \label{sayfa} 
$$
\begin{array}{lll}
F_{(1,0)(2)}, & F_{(2,0)(1)}, & F_{(0)(2,1)}, \\ 
F_{(0)(2)}, & F_{(1)(2)}, & F_{(0)(1)}. 
\end{array}
$$
For all $x_1\in NG_1,~y_2\in NG_2,$ the corresponding generators of $N_3$ are:\label{sayfa1} 
$$
\begin{array}{lll}
F_{(1,0)(2)}(x_1, y_2) & = & [s_1s_0x_1, s_2y_2]{~} [s_2y_2, s_2s_{0}x_1] \\ 
F_{(2,0)(1)}(x_1, y_2) & = & [s_2s_0x_1, s_1y_2]{~} [s_1y_2, s_2s_1x_1]{~}  
[s_2s_1x_1, s_2y_2]{~} [s_2y_2, s_2s_0x_1]  
\end{array}
$$
and all $x_2\in NG_2,~y_1\in NG_1,$
$$
\begin{array}{lll}
F_{(0)(2,1)}(x_2, y_1) & = & [s_0x_2, s_2s_1y_1]{~}[s_2s_1y_1, s_1x_2]{~}
[s_2x_2, s_2s_1y_1] 
\end{array}
$$
whilst for all $x_2 , y_2 \in NG_2$,\label{sayfa1b} 
$$
\begin{array}{lll}
F_{(0)(1)}(x_2 , y_2) & = & [s_0x_2 , s_1y_2]{~} [s_1y_2 , s_1x_2]{~} [s_2x_2 ,s_2y_2] \\ 
F_{(0)(2)}(x_2 , y_2) & = & [s_0x_2 , s_2y_2]  \\ 
F_{(1)(2)}(x_2 , y_2) & = & [s_1x_2 , s_2y_2]{~}[s_2y_2 , s_2x_2]. 
\end{array}
$$
Our aim in this paper is to prove that the images of these pairings generate 
$NG_n\cap D_n.$ More precisely:
\begin{thm}~~{\bf ( Theorem A )}\label{10}
Let {\bf G} be a simplicial group and for $n>1,$ let $D_n$ the 
subgroup in $G_n$ generated by degenerate elements.  
Let $N_n^G$ be the normal subgroup generated by elements of the form 
$$
F_{\alpha ,\beta }(x_\alpha , y_\beta )\qquad \text{with }(\alpha
,\beta )\in P(n) 
$$
where $x_{\alpha} \in NG_{n-\#\alpha },\ y_{\beta} \in NG_{n-\#\beta }$.
Then 
$$
NG_n\cap D_n = N_n^G\cap D_n. 
$$
\end{thm}
As a corollary we, of course, have that the image of $N_n^G\cap D_n$ is equal to the image of $NG_n\cap D_n$ 
i.e., $\partial _n(N_n\cap D_n)=\partial _n(NG_n\cap D_n)$.\\

The proof of 2.1 is given in the next section after some preparatory lemmas. Here we
restrict to the case $n=2$ by way of illustration. In their paper \cite{bl1}, 
Brown and Loday proved a lemma:

\begin{lem}~\cite{bl1}
Let ${\bf G}$ be a simplicial group such that $G_2 = D_2$ is generated by 
degenerate elements. Then in the Moore complex of ~${\bf NG}$ we have 
$\partial_2NG_2 = \partial_2N_2$ where $N_2$ is the normal subgroup of $G_2$ generated
by elements of the form
$$
\begin{array}{lll}
F_{(0)(1)}(x_1,y_1) &=&  [s_0x_1,s_1y_1]{~}[s_1y_1, s_1x_1]\\
\end{array}
$$
with $x_1,y_1\in NG_1.$
\end{lem}
This is, of course, a trivial consequence of Theorem A and their proof
inspired that of the more general theorem given here.

\textbf{Remark:}

AN unknown referee made the interesting observation that if $x$ is a Moore
cycle, so $\partial x = 0$, then $F_{(0)(1)}(x,x)$ is one also.  Thus
$F_{(0)(1)}$ induces an operation $\pi_\ast(G) \rightarrow \pi_{\ast +1}(G)$.
Geometrically this operation can be described as the $\eta$-operation given by 
the composition $S^{m+1}\stackrel{\eta}'{\rightarrow}S^m
\stackrel{x}{\rightarrow}G$, where $\eta$ is the (suspension of) the Hopf map.
 The geometric interpretation of the $F_{\alpha \beta(1)}$ in general would
 seem to be quite important but the authors have as yet little idea what it
 might be.
\section{Elements of $N_n$ and properties of the pairings}
In the following we analyse various types of elements in $N_n$ and show that
products of them give elements that we want in giving an alternative description of 
$NG_n$.
\begin{lem}\label{7}
Given $x_\alpha \in NG_{n-\#\alpha },~y_\beta \in NG_{n-\#\beta }$ with $%
\alpha =(i_l,\ldots ,i_1),$ $\beta =(j_m,\ldots ,j_1)\in S(n).$ If $\alpha
\cap \beta =\emptyset $ with $\beta <\alpha $ and $v=[s_\alpha x_\alpha,
s_\beta y_\beta],$ then

(i) \ \ if $k\leq i_1,$ then $p_k(v)=v,$ 

(ii)\ \ \thinspace if $k>i_l+1\,$or $k>j_m+1,$ then $p_k(v)=v,$ 

(iii)\ \ if $k\in \{j_1,\ldots ,j_m \}$ and $k=i_r+1$ for some $r,$ then for \\
${\alpha'} = (i_l,\ldots,i_r+1,i_r-1,\ldots, i_1)$ and
${\beta} = (j_m,j_{m-1},\ldots, j_1),$
$$
p_k(v)=[s_\alpha x_\alpha, s_\beta y_\beta]{~}[s_{\alpha^{\prime}}
x_{\alpha}, s_{\beta}y_{\beta}]^{-1},
$$

(iv)\ \ \ if $k\in \{i_1, \ldots, i_l \}$ and $k=j_s+1$ for some $s,$ then 
for \\
${\beta'} = (j_m,\ldots,j_s+1,j_s-1\ldots, j_1)$
$$
\begin{array}{lll}
p_k(v) & = & [s_\alpha x_\alpha, s_\beta y_\beta ]{~}
[s_{\alpha}x_{\alpha}, s_{\beta^{\prime}}y_{\beta}]^{-1} \\ 
& = & vv^{\prime}  
\end{array}
$$
where $v^{'}\in G_{n-1}$ and $0\leq k\leq n-1,$ 

(v)\ \ \ if $k = j_m+1$  (or $k = i_l+1$) then 
$$
\begin{array}{lll}
p_k(v) & = & vs_k(v_k)^{-1}\\
& = & [s_\alpha x_\alpha, s_\beta y_\beta]s_k(v_k)^{-1}
\end{array}
$$
where $s_k(v_k)^{-1} = [s_{\beta^{\prime}}y_{\beta}, s_{\alpha}x_{\alpha}]$
$(\text{or}\quad  s_k(v_k)^{-1} = [s_{\beta}y_{\beta}, s_{\alpha^{\prime}}x_{\alpha}])$
with respect to $k= j_m+1$ (and $k= i_l+1$ respectively) and for new strings 
${\alpha'}$ and ${\beta'},$ 

(vi)\ \ \ if $k = j_1+1,$ then 
$$
\begin{array}{lll}
p_k(v) & = & vs_k(v_k)\\
& = & [s_\alpha x_\alpha, s_\beta y_\beta]s_k(v_k)^{-1}
\end{array}
$$
where $s_k(v_k)^{-1}=[s_{\beta^{\prime}}y_{\beta}, s_{\alpha}x_{\alpha}]$

(vii)\ \ \ if $k \in \{j_1,\ldots,j_m,j_{m+1}\}$ and $k=i_t+1$ for some $t,$ then
$$
\begin{array}{lll}
p_k(v) & = & [s_\alpha x_\alpha, s_\beta y_\beta]
[s_{\beta'}y_\beta, s_{\alpha'}x_\alpha]
\end{array}
$$
where $0\leq k \leq n-1.$
\end{lem}
\begin{pf}
Assume $\beta <\alpha $ and $\alpha \cap \beta =\emptyset $ which implies $%
i_1<j_1.$ In the range $0\leq k\leq i_1,$ 
$$
\begin{array}{lll}
p_k(v) & = & [s_\alpha x_\alpha, s_\beta y_\beta ]{~}[s_kd_ks_{i_{m}}\ldots 
s_{i_{1}}x, s_kd_ks_{j_{m}}\ldots s_{j_{1}}y]^{-1} \\
& = & [s_\alpha (x_\alpha ), s_\beta (y_\beta )]{~}[s_{j_{m}-1}\ldots s_{j_{1}-1}s_kd_ky,
s_{i_{m}-1} \ldots s_{i_{1}-1}s_kd_kx] \\
& = & [s_\alpha x_\alpha, s_\beta y_\beta]\qquad \text{since }d_k(x_\alpha
)=1 \ \ or \ \ d_k(y_\beta)=1. \\
& = & v .\\
\end{array}
$$
Similarly if $k>i_l+1,$ then%
$$
\begin{array}{lll}
p_k(v) & = & [s_\alpha x_\alpha, s_\beta y_\beta ]{~}[s_{j_{m}}\ldots s_{j_{1}}
s_{k-m}d_{k-m}y, s_{i_{m}}\ldots s_{i_{1}}s_{k-m}d_{k-m}x] \\
& = & [s_\alpha x_\alpha, s_\beta y_\beta ]\qquad \text{since }%
d_{k-m}(y_\beta )=1 \ \ or \ \ d_{k-l}(x_\alpha)=1.\\ 
& = & v .
\end{array}
$$
Clearly the same sort of argument works if $k>j_{m+1}.$
\noindent If $k\in \{j_1,\ldots ,j_m,j_{m+1} \}$ and $k=j_t+1$ for some $t,$
then%
$$
\begin{array}{lll}
p_k(v) & = & [s_\alpha x_\alpha, s_\beta y_\beta ]{~}[s_kd_ks_\beta
y_\beta , s_kd_ks_\alpha x_\alpha] \\  
& = & [s_\alpha x_\alpha, s_\beta y_\beta]{~}[s_{\beta^{\prime}}
x_{\beta^{\prime}}, s_{\alpha^{\prime}}y_{\alpha^{\prime}}]\\ 
\end{array}
$$
\hfill\end{pf}
\begin{lem}\label{8}
If $\alpha \cap \beta =\emptyset $ and $ \beta<\alpha ,$ then 
$$
p_{l}\ldots p_1{~} [s_\alpha x_\alpha , s_\beta y_\beta]= [s_\alpha x_\alpha
, s_\beta y_\beta ] \prod\limits_{i=1}^{l}s_i(z_i)^{-1} 
$$
where $z_i \in \bigcap\limits_{j=0}^{i-1}$\text{\rm Ker}$d_j\subset G_{n-1}$ and $l\in [n-1].$
\end{lem}
\begin{pf}
By induction on $l.$
\end{pf}
\begin{lem}\label{9}
Let $x_\alpha \in NG_{n-\#\alpha },~y_\beta \in NG_{n-\#\beta }$ with $%
\alpha ,\beta \in S(n),$ then 
$$
s_\alpha x_\alpha s_\beta y_\beta s_\alpha (x_\alpha )^{-1}
=s_{\alpha \cap \beta }z_{\alpha \cap
\beta } 
$$
where $z_{\alpha \cap \beta }$ has the form $s_{\bar{\alpha}}x_\alpha
s_{\bar{\beta}}y_\beta s_{\bar{\alpha}}(x_\alpha)^{-1}
$ and $\bar{\alpha} \cap \bar{\beta}=\emptyset .$
\end{lem}
\begin{pf}
If $\alpha \cap \beta =\emptyset ,$ then this is trivially true. Assume $%
\#(\alpha \cap \beta )=t,$ with $t\in {\mathbb N}.$ Take $\alpha =(i_l,\ldots ,i_1)$ ~and ~$%
\beta =(j_m,\ldots ,j_1)$ with $\alpha \cap \beta =(k_t,\ldots ,k_1),$ 
$$
s_\alpha x_\alpha = s_{i_l}\ldots s_{k_t}\ldots s_{i_1}x_\alpha ~~~\text{%
and}~~~s_\beta y_\beta = s_{j_m}\ldots s_{k_t}\ldots s_{j_1}y_\beta . 
$$
Using repeatedly the simplicial axiom $s_as_b=s_bs_{a-1}$ for $b<a$ until
 $s_{k_t}\ldots s_{k_1}$ is at the beginning of the string, one
gets the following 
$$
s_\alpha x_\alpha = s_{k_t\ldots k_1}(s_{\bar{\alpha} }x_\alpha )~~~ 
\text{and}~~~s_\beta y_\beta = s_{k_t\ldots k_1}(s_{\bar{\beta}}y_\beta ). 
$$
Multiplying these expressions together gives 
$$
\begin{array}{lll}
s_\alpha x_\alpha s_\beta y_\beta s_\alpha (x_\alpha )^{-1}
& = & s_{k_t}\ldots s_{k_1}(s_{\bar{\alpha}}
x_\alpha )s_{k_t}\ldots s_{k_1}(s_{\bar{\beta}}y_\beta ) 
s_{k_t}\ldots s_{k_1}(s_{\bar{\alpha}}(x_\alpha)^{-1}) \\ 
& = & s_{k_t\ldots k_1}(s_{\bar{\alpha}}x_\alpha s_{\bar{\beta}} 
y_\beta s_{\bar{\alpha}}(x_\alpha)^{-1})\\  
& = & s_{\alpha \cap \beta }(z_{\alpha \cap \beta }), 
\end{array}
$$
where $z_{\alpha \cap \beta }= s_{\bar{\alpha}}x_\alpha s_{\bar{\beta}}
y_\beta s_{\bar{\alpha}}(x_\alpha)^{-1}
\in {NG_{n-\#(\alpha \cap \beta )}}$ and where 
$\bar{\alpha} = (i_l-t,\ldots, k_t+1-t,\dots,i_1)$ and 
$\bar{\beta} =(j_m-t,\ldots, {k_t'}+1-t,\dots,j_1).$ 
Hence $\bar{\alpha} \cap \bar{\beta} 
=\emptyset_{n-\#(\alpha \cap \beta )}.$ 
Moreover $\bar{\alpha} <\alpha $ and $~\bar{\beta} 
<\beta $ as $\#\bar{\alpha}<\#\alpha $ and $~\#\bar{\beta} 
<\#\beta.$
\end{pf}

Suppose $\alpha= (i_s,\ldots,i_1)\in S(m)$ and $\gamma:[n]\longrightarrow
[m]\in S(n,m)$. Define $\gamma_{\ast}(\alpha)$   by 
$s_{\gamma_{\ast}(\alpha)}= s_{\gamma}s_{\alpha}$. 
\begin{cor}\label{baba}
Let  $\beta \leq \alpha$ and 
$\gamma:[n]\longrightarrow [m].$ 
Then $\gamma_{{\ast}(\beta)}\leq\gamma_{{\ast}(\alpha)}\Longleftrightarrow\beta\leq\alpha,$
where $\gamma_{{\ast}(\alpha)},\ \gamma_{{\ast}(\beta)}\in S(n).$ $\square$
\end{cor}
The following lemma is proved similarly.
\begin{lem}\label{anne}
For $m\leq n,$ suppose given in $G_m$ an element 
$$
g = \prod\limits_{{\beta'}\leq{\gamma'}\leq{\alpha'}}s_{\gamma'}(z_{\gamma'})
$$
and $s_{\delta}: G_m\longrightarrow G_n.$ Then setting $\alpha, \ \beta\in S(n)$ 
such that\\
$s_{\delta}s_{\alpha'} = s_{\alpha}, \ s_{\delta}s_{\beta'} = s_{\beta}$
$$
s_{\delta}(g) = \prod\limits_{{\beta}\leq{\gamma}\leq{\alpha}}s_{\gamma}(z_{\gamma})
$$
for some elements $z_{\gamma}\in NG_{n-\#{\gamma}}$ and where 
$s_{\delta}s_{\gamma'} = s_{\gamma}.$ $\square$
\end{lem}
{\noindent \bf Proof of Theorem A :}\\
From Proposition \ref{2}, $G_n$ is isomorphic to%
$$
NG_n\rtimes s_{n-1}NG_{n-1}\rtimes s_{n-2}NG_{n-1}\rtimes \ldots \rtimes  
s_{n-1}s_{n-2}\dots s_0NG_0. 
$$
Similarly
$D_n$ is isomorphic to\\
$(NG_n\cap D_n)\rtimes s_{n-1}NG_{n-1}\rtimes s_{n-2}NG_{n-1}\rtimes \ldots \rtimes
s_{n-1}s_{n-2}\dots s_0NG_0.$~~
Hence any element $g$ in $D_n$ can be written in the following form%
$$
g = g_ns_{n-1}(y_{n-1})s_{n-2}(y_{n-1}^{\prime
})s_{n-1}s_{n-2}(y_{n-2})\ldots s_{n-1}s_{n-2}\ldots s_0(y_0),\text{\quad} 
$$
with $g_n\in NG_n\cap D_n,$ $y_{n-1},y_{n-1}^{\prime }$ $\in NG_{n-1},$ $y_{n-2}\in
NG_{n-2},$ $y_0$ $\in NG_0$ etc.\\
To simplify the notation a little, we will assume that $G_n = D_n$, so that
$N_n \subset D_n$.  The general case would replace $NG_n$ by $NG_n\cap D_n$ and 
similarly $N_n$ by $N_n\cap D_n$ from here on.

As it is easily checked that $N_n\subseteq NG_n\cap D_n,$ 
it is enough to prove that any element in $D_n/N_n$ can be 
written in the form
$$
s_{n-1}(y_{n-1})s_{n-2}(y_{n-1}^{\prime })s_{n-1}s_{n-2}(y_{n-2})\ldots
s_{n-1}s_{n-2}\ldots s_0(y_0)N_n 
$$
that is, for any $g\in D_n$,
$$
gN_n=s_{n-1}(y_{n-1})s_{n-2}(y_{n-1}^{\prime })\ldots
s_{n-1}s_{n-2}\ldots s_0(y_0)N_n. 
$$
for some $y_{n-1}\in NG_{n-1},$ etc. We refer to this as the standard form of $gN_n.$ \newline
If $g\in D_n,$ it is a product of  degeneracies.
If $g$ is itself a degenerate element, it is obvious 
that it is a product of elements in the semidirect factors, 
$s_\beta (G_{n-\#\beta }),$ ~$\beta\in S(n)- \{\emptyset_n\}.$
\newline
Assume therefore that
provided an element $g$ can be written as a product of $k-1$ degeneracies of 
this form, then it has the desired form  modulo $N_n.$ Now for an element $g$ which
needs $k$ degenerate elements, we have 
$$
g=s_{\alpha}x_\alpha g^{\prime }
\qquad \text{with }x_\alpha \in NG_{n-\#\alpha} 
$$
where $g^{\prime }$ needs fewer than $k$ and so%
$$
\begin{array}{lll}
gN_n & = & s_\alpha x_\alpha g^{\prime }N_n \\  
& = & s_\alpha x_\alpha (s_{n-1}(y_{n-1})s_{n-2}(y_{n-1}^{\prime })\ldots
s_{n-1}s_{n-2}\ldots s_0(y_0))N_n. 
\end{array}
$$
We prove that this can be rewritten in the desired form  $\text{mod}~N_n$
by using induction on $\alpha$ within the linearly ordered set 
$S(n)- \{\emptyset_n\}.$ \\
If $\alpha =(n-1)$, then 
$$
gN_n = s_{n-1}(xy_{n-1})s_{n-2}(y_{n-1}^{\prime })\ldots
s_{n-1}s_{n-2}\ldots s_0(y_0) N_n 
$$
where $x\in NG_{n-1}$ and $(xy_{n-1})\in NG_{n-1}.$ \\
If $\alpha =(n-2),$ then since 
$$
F_{(n-2)(n-1)}(x_{n-1}, y_{n-1}) = [s_{n-2}x_{n-1},~ s_{n-1}y_{n-1}]
~[s_{n-1}y_{n-1},~ s_{n-1}x_{n-1}]
$$
and
$$
s_{n-2}(x_{n-1})s_{n-1}(y_{n-1})s_{n-2}(x_{n-1})^{-1} \equiv s_{n-1}(x_{n-1}y_{n-1}x^{-1}_{n-1})
\quad\text{mod}~N_n,
$$
we have
$$
\begin{array}{lll}
gN_n & =& (s_{n-2}(x_{n-1})s_{n-1}(y_{n-1})s_{n-2}(x_{n-1})^{-1})
s_{n-2}(x_{n-1})s_{n-2}(y_{n-1}^{\prime }) \\
& & \ldots s_{n-1}s_{n-2}\ldots s_0(y_0) N_n \\
& = & s_{n-1}(xyx^{-1}_{n-1})s_{n-2}(xy_{n-1}^{\prime })\ldots
s_{n-1}s_{n-2}\ldots s_0(y_0) N_n
\end{array}
$$
where $x_{n-1}, y_{n-1}\in NG_{n-1}$ so $(xyx_{n-1}^{-1}), ~~ 
(xy_{n-1}^{\prime })\in NG_{n-1}.$ \\
In general we need to sort $s_{\alpha}x_{\alpha}$ into its correct 
place in the product but in so doing will conjugate earlier terms in the product 
as happened in the case  $\alpha= (n-2)$ above. Each of these terms must be 
shown to consist only of subterms of types we have already dealt with, that is 
further to the left in the standard form of the product. Explicitly we assume 
that we can do this sorting for any term $s_{\gamma}x_{\gamma}$ with 
$\gamma < \alpha$ and examine
$$ 
\begin{array}{lll}
gN_n & =& s_{\alpha}x_{\alpha}(s_{n-1}(y_{n-1})s_{n-2}(y_{n-1}^{\prime })\ldots
s_{n-1}s_{n-2}\ldots s_0(y_0))N_n \\
& = &s_{\alpha}x_{\alpha} 
\prod\limits_{\beta\in S(n)-\{\emptyset_{n}\}}s_{\beta}y_{\beta}N_n \\
gN_n & = &\prod\limits_{\alpha>\beta}s_\alpha x_\alpha s_\beta y_{\beta}s_\alpha (x_\alpha )^{-1}
\cdot s_\beta (xy)_\beta \cdot~\prod\limits_{\beta >\alpha}s_\beta y_{\beta}N_n
\end{array}
$$
where $\beta\in S(n)-{\emptyset_n}$ and $\alpha>\beta$ with respect to 
the lexicographic ordering in $S(n).$ \\
We  look at products of the following type 
$$
s_\alpha x_\alpha s_\beta y_\beta 
s_\alpha (x_\alpha)^{-1}\qquad (\ast) 
$$
and we want to show that these can always be written in the form
$$
\prod\limits_{\gamma\leq\beta}s_{\gamma}(z_{\gamma})
$$
for some $z_{\gamma}\in NG_{n-\#\gamma}.$ This will mean that 
we already know how to sort all the terms that arise since none occur
`to the right of' $\beta$ in the lexicographic order in the product.\\
We check this product case by case as follows: \\
If $\alpha\cap\beta = \emptyset,$ then by Lemma~\ref{8}, 
$$
s_\alpha x_\alpha s_\beta y_\beta s_\alpha (x_\alpha )^{-1}
\equiv\prod\limits_{k=l}^{i_1+1}s_k(z_k)s_\beta y_\beta \quad\text{mod}~N_n,
$$
where $\beta\in S(n)-\{\emptyset_n\}.$
Now we need to show that each $s_k(z_k)$ is made up  of terms $s_\mu(z_\mu)$
with $(z_\mu)\in NG_{n-\#\mu},$ $\mu\leq\beta.$ (We will use the notation of Lemma ~\ref{7}.) 
Since $\alpha>\beta$ then $i_1\leq j_1.$ We have 
$$
z_k = \prod\limits_{\mu\leq(k-1)}s_\mu(z_\mu)
\quad\text{since}\quad z_k\in \bigcap\limits_{j=0}^{k-1}\text{Ker}d_j
$$
so
$$
s_k(z_k) = \prod\limits_{\mu\leq(k-1)\leq(i_1)}s_ks_\mu(z_\mu)
$$
where we write $\mu =(m_1,\ldots,m_r)$ so we have $k-1\leq m_1$ or $k\leq m_1+1.$ \\
We compare $k$ with $m_1$ and $m_2$: either \\
(a) ~$k = m_1+1 < m_2;$ \\
(b) ~$k = m_1+1 = m_2 ;$ \\
$\spreaddiagramrows{-1.2pc} \spreaddiagramcolumns{-1.2pc}
\def\objectstyle{\ssize} \def\labelstyle{\ssize}
(c) ~k=m_1 = \left\{\begin{array}{ll}
\text{and}& m_2 =m_1+1, \text{or}\\
\text{and}& m_2 > m_1+1,
\end{array}\right.$\\ 
or \\
(d) $k< m_1.$ \\
Thus \\
\begin{equation*}
s_ks_{\mu}= s_{m_r+1}\ldots s_ks_{m_2}s_{m_1} =\left\{\begin{array}{lll}
s_{m_r+1}\ldots s_{m_2+1}s_ks_{m_1}& \quad\text{cases (a) and (b),}\\
s_{m_r+1}\ldots s_{m_1+1}s_{m_1} &\quad\text{case (c),}\\
s_{m_r+1}\ldots s_{m_1+1}s_{k} &\quad\text{case (d),}
\end{array}\right.
\end{equation*}
so in each case $s_ks_{\mu} = s_{\mu'}$ where ${m_1'} = \text{min}\{k,m_1\}.$
We compare ${\mu'}$ with $\beta.$ If ${m_1'} > j_1,$ then ${\mu'}\leq \beta.$
If ${m_1'} =j_1,$ then either $k=j_1$ or $m_1 =j_1$ then $k = j_1+1.$ Thus we need 
to show 
$$
w=s_{j+1}d_{j+1}[s_{\beta} y_{\beta} ,~s_{\alpha} x_{\alpha} ] = 
\prod\limits_{\vartheta\leq\beta}s_{\vartheta}(z_{\vartheta}).
$$ 
There are two cases:\\
(i){~}If $j_1\in \alpha,$ then 
$$
w=[s_{\beta'} y_{\beta} ,~s_{\alpha} x_{\alpha} ] = 
\prod\limits_{\vartheta\leq{\beta'}\leq\beta}s_{\vartheta}(z_{\vartheta})
$$
where $\beta\leq{\beta'}$ and since $j_1+1 =j_s\notin\beta,$ then 
${\beta'} = \{j_{m},\ldots, j_{s+1}+1, j_s, j_{s-1}, \ldots,j_1\},$ \\
(ii){~}If $j_1+1\in \beta,$ then 
$$
w=[s_{\beta} y_{\beta} ,~s_{\alpha'} x_{\alpha} ] = 
\prod\limits_{\vartheta\leq\beta}s_{\vartheta}(z_{\vartheta})
$$
where $\alpha\leq{\alpha'}$and since $j_1+1= i_r\notin\alpha,$ then 
${\alpha'} = \{i_{l}, \ldots, i_{r+1}+1, i_r, i_{r-1}, \ldots,i_1\}.$ \\ 
Both cases are covered by the induction hypothesis. Both cases can thus be written 
$$
\prod\limits_{\vartheta\leq\beta}s_{\vartheta}(z_{\vartheta}).
$$
We have $z_{\vartheta}\in NG_{n-\#\vartheta}$ and for 
some $r$ and $s$ then it can be written for above cases
$$
s_k(z_k) = \prod\limits_{{\emptyset}_n\leq{\gamma'}\leq(k)\leq(i_1)}s_{\gamma'}(z_{\gamma})
$$
where $s_ks_\gamma = s_{\gamma'},$  ${\gamma'}\leq\beta$ and $z_{\gamma'}\in NG_{n-\#\gamma'}$
so 
$$
gN_n = \prod\limits_{{\gamma'}\leq\beta}s_{\gamma'}(z_{\gamma'})\cdot s_{\alpha}(xy_{\alpha})\cdot
\prod\limits_{{\alpha}<\beta}s_{\beta} y_{\beta}N_n 
$$
as required.\\
If $\alpha \cap \beta \neq \emptyset $, then one gets, from Lemma \ref{9},
the following 
$$
s_\alpha  x_\alpha  s_\beta  y_\beta  s_\alpha (x_\alpha)^{-1}
= s_{\alpha \cap \beta }(s_{\bar{\alpha}} x_\alpha  s_{\bar{\beta}} y_\beta  
s_{\bar{\alpha}}(x_\alpha )^{-1}) 
$$
where $\bar{\alpha}> \bar{\beta},$  $\bar{\alpha}\cap\bar{\beta}\in \emptyset_{n-\#(\alpha \cap \beta)}.$
Using Lemma ~\ref{anne} in dimension $n-\#(\alpha \cap \beta)$ and Corollary ~\ref{baba} then we have
$$
(s_{\bar{\alpha}} x_\alpha  s_{\bar{\beta}} y_\beta  s_{\bar{\alpha}}(x_\alpha )^{-1})\equiv 
\prod\limits_{\theta\in[\emptyset_m, ~\bar{\beta}]}s_{\theta}(z_{\theta})
$$
hence 
$$
\begin{array}{lll}
s_\alpha  x_\alpha  s_\beta  y_\beta  s_\alpha (x_\alpha)^{-1} & = &\prod\limits_{\theta\in 
[\emptyset_m, ~\bar{\beta}]}
s_{\alpha \cap \beta }s_{\theta}(z_{\theta}) \\
&=& \prod\limits_{\eta\in[\alpha\cap\beta, ~\beta]}s_{\eta}(z_{\eta})
\end{array}
$$
where $z_{\eta} \in NG_{n-\#\eta},$ 
\begin{displaymath}
z_{\nu} =\left\{\begin{array}{ll}
z_{\eta} & \text{if $\eta = \gamma_{\ast}(\eta)$} \\
1 & \text{otherwise}
\end{array}\right.
\end{displaymath} 
and $s_{\alpha \cap \beta }s_{\theta} = s_\eta.$ Then $gN_n$ can be written
$$
\begin{array}{lll}
gN_n & = &\prod\limits_{\nu\in[\alpha \cap \beta, ~\beta]}s_{\nu}(z_{\nu})
\cdot
s_\beta  (xy)_\beta \cdot\prod\limits_{\beta>\alpha}s_\beta  y_\beta  N_n.
\end{array}
$$
where $s_\eta = s_\nu.$
Thus we have shown that every product can be rewritten in the required
form modulo $N_n$, so in general, $N_n\cap D_n = NG_n\cap D_n.$  $\square$
\section{Applications and implications}
Kan introduced the notation of a $CW$-basis for a free simplicial group and used this to proved that 
free simplicial groups model all connected homotopy types. The idea is that as one adds cells to the 
$CW$-complex one adds new generators to the free simplicial group, but one does this within the 
Moore complex so the new simplices have all but their last face at the identity element in the next 
dimension down. In homotopy types where there are few such non degenerate generators or where these 
generators are `generated' in a simple way then the methods behind Theorem A raise the hope of finding 
a detailed presentation of the segments of the homotopy type between those dimensions in which there 
are non-degenerate generators. The means of presenting this information may vary with the context, but 
one set of fairly compact methods comes from the crossed algebraic techniques pioneered by J. H. C. 
Whitehead in \cite{whitehead}. (Modern references for this and for more recent developments can 
conveniently be found in the survey article by Baues \cite{baues3}.)
\subsection{Crossed complexes}
As an illustration we examine the impact of Theorem A on the links between simplicial groups and the 
homotopy systems of Whitehead, more exactly the connected crossed complexes of Brown and Higgins 
(cf. \cite{bh1} and \cite{bh2}) or the crossed chain complexes of Baues (cf. \cite{baues2} and 
\cite{baues3}) as no freeness assumptions will be made here.

Let $G$ be a group, then a $G$-group is a group $H$ together with a given action of $G$ on $H,$ that is a
homomorphism from $G$ to $\text{Aut}(H).$ 
\begin{defn}(cf. Baues, \cite{baues3} p.22)
A crossed complex  $\rho$ is a sequence
$$
\diagram
\rto^{d_4}&\rho_3\rto^{d_3}&\rho_2\rto^{d_2}&\rho_1
\enddiagram
$$
of homomorphisms between $\rho_1$-groups where $d_2$ is a crossed module and $\rho_n,$ $n\ge 3$ is 
abelian and a $\pi_1$-module via the action of  $\rho_1,$ where $\pi_1=\text{cokernel}(d_2).$ Moreover 
$d_{n-1}d_n = 0$ for $n\ge 3.$
\end{defn}

It is known (cf. Ashley \cite{ashley}, Carrasco and Cegerra \cite{c:c} or Ehlers and Porter \cite{ep} 
and the references therein) that crossed complexes correspond, via a nerve-type functor, to simplicial 
groups with a `thin' structure. Each simplicial group is a Kan complex as a well known algorithm gives 
a filler for any horn. A Kan complex, $K,$ is a $T$-complex if there is for each $n$ a subset $T_n$ of 
$K_n,$ made up of so called `thin' elements, such that any horn has a unique thin filler and two other 
more technical conditions hold (cf. Ashley \cite{ashley}). A simplicial $T$-complex which is also a 
simplicial group is a group $T$-complex provided in each dimension $T_n$ is a subgroup of $K_n.$ In 
this case one easily checks that $T_n$ must be $D_n$ the subgroup of $K_n$ generated by the degenerate 
elements. 
\begin{prop}\cite{ashley}
A simplicial group ${\bf G}$ has ${\bf NG}$ a crossed complex if and only if for each $n\geq 1,$ 
$NG_n\cap D_n$ is trivial.~$\square$
\end{prop}

The idea of the proof is that two $D_n$ fillers for the same horn must differ by an element of 
$NG_n\cap D_n,$ so uniqueness corresponds to the simplicial group being a group $T$-complex. 
The final part uses Ashley's  equivalence between group $T$-complexes and crossed complexes.

Carrasco and Cegarra used this in \cite{c:c} to prove (p. 215) that a simplicial group has Moore 
complex a crossed complex if and only if their pairings vanish (these are similarly defined to 
those used  here but are based on products rather than commutators). Similarly we have:
\begin{cor}
A necessary and sufficient condition that a simplicial group ${\bf G}$ has ${\bf NG}$ a crossed 
complex is that for all $n$ and all $\alpha,\beta\in P(n),\ F_{\alpha, \beta}^{n}(x,y)$ is trivial 
for all pairs $(x,y).$
\end{cor}
\begin{pf}
Since the $F_{\alpha, \beta}^{n}(x,y)$ normally generate $NG_n\cap D_n,$ this is immediate.
\end{pf}

The importance of this result is probably for the interpretation of the  $F_{\alpha, \beta}^{n}$ 
as their vanishing has a great simplifying   effect on the Moore complex.
\subsection{$\sum K(\pi,1)$}
As mentioned earlier Brown and Loday used their generalised van Kampen theorem to calculate 
$\pi_3\sum K(\pi, 1)$ as $\text{Ker}(\pi\otimes\pi\longrightarrow\pi),$ the kernel of the commutator 
map. Jie Wu (\cite{wu} Theorem 5.9) proves that for any group, $\pi,$ and set of generators 
$\{x^{\alpha} \mid \alpha\in J\}$ for $\pi,$ then for $n\ne 1,$ $\pi_{n+2}(\sum K(\pi, 1))$ is 
isomorphic to the center  of the quotient group of the free product
$$
\coprod\limits_{0\leq j\leq n}(\pi)
$$
modulo the relation
$$
[y_{i_1}^{(\alpha_1)\varepsilon_1}, y_{i_2}^{(\alpha_2)\varepsilon_2},
\ldots, y_{i_t}^{(\alpha_t)\varepsilon_t}]
$$
where $\{i_1,\ldots,i_t\}=\{ -1,0,\ldots, n\}$ as sets, where $(\pi)_j$ is a copy of $\pi$ 
with generators $\{x_{j}^{(\alpha)}\mid \alpha\in J\}, \ \varepsilon_j= \pm 1,$ 
$y_{-1}^{(\alpha)}= x_0^{(\alpha)^{-1}}, \ y_{j}^{(\alpha)}= x_j^{(\alpha)}x_{j+1}^{(\alpha)^{-1}}$ 
for $1\leq i\leq {n-1}$ and $y_n^{(\alpha)}= x_{n}^{(\alpha)},$ and finally the commutator bracket 
$[\ldots]$ runs over all the commutator bracket arrangements of weight $t$ for each $t.$

Wu's methods rely on using a construction he attributes to Carlsson \cite{carlsson}. This gives a 
simplicial group $F^{\pi}(S^1)$ that has $\pi_{n+2}\sum K(\pi,1)\cong \Omega\sum K(\pi,1)\cong 
\pi_{n+1}F^{G}(S^1).$ Our Theorem A above provides a link between Wu's methods and the Brown-Loday 
result. We will explore this link to some extent but cannot as yet retrieve the Brown-Loday result by
purely algebraic methods. Potentially  however this might yield a tensor-like description in dimension 
$4$ and higher, but we will not explore that here.

Although Carlsson introduced the construction $F^G(X)$ in 1984, the construction is essentially the 
same as the tensoring operation used by Quillen and others. Working in the
simplicially enriched category of 
simplicial groups, there is a tensor operation defined as follows: let $K$ be a simplicial set and 
$G_1, \ G_2$ simplicial groups. The simplicial group $G_1\bar{\otimes}K$ has the universal property 
given by the natural isomorphism
$$
S(K, SGp(G_1, G_2))\cong SGp(G_1\bar{\otimes}K, G_2).
$$
The category of simplicial groups is also enriched over $S_{\ast},$ the category of pointed 
simplicial sets. We define $G_1\bar{\wedge}K$ by
$$
S_{\ast}(K, SGp(G_1, G_2))\cong SGp(G_1\bar{\wedge}K, G_2).
$$
There is an isomorphism $F^{G}(K)\cong G\bar{\wedge}K.$ The advantage of this
approach is that  it makes it clear that $\bar{\wedge}$ generalises $\wedge$
just as $\bar\otimes$ generalises $\times$
\begin{lem}
If $f: G\longrightarrow H$ is a morphism of simplicial groups, it induces 
$f\bar{\wedge}K : G\bar{\wedge} K\longrightarrow H\bar{\wedge} K,$ moreover if $f$ is a weak 
homotopy equivalence, so is  $f\bar{\wedge}K.$
\end{lem}
\begin{pf}
As Carlsson noted, $(G\bar{\wedge}K)_n$ is
$$
\coprod\limits_{x\in K_n}(G_n)_{x}/(G_n)_{\ast}
$$
and is thus the diagonal of a bisimplicial group having
$\coprod\limits\{(G_m)_{x}\mid x\in K_n\setminus\{{\ast}\}\}$ in its $(m,n)$-position.
A simple spectral sequence argument, or direct manipulation, completes the proof.
\end{pf}
\begin{prop}
There is a natural weak homotopy equivalence
$$
\Omega\sum K(\pi,1)\simeq K(\pi,0)\bar{\wedge}S^{1}
$$
where $K(\pi,0)$ is the constant simplicial group with value, $\pi,$ $S^1$ is the simplicial 
$1$-sphere and $\sum$ is reduced suspension.
\end{prop}
\begin{pf}
As Kan's loop group functor models $\Omega$ and $K\wedge S^1$ the suspension,
$$
\Omega(\sum K(\pi, 1))\simeq G(K(\pi,1)\wedge S^1)
$$
then the adjunction between $G$ and the classifying space functor $\bar{W}$
gives for $K$, $ L$, arbitrary pointed simplicial sets, and ${\bf H}$  an arbitrary 
simplicial group, the natural isomorphisms 
$$
\begin{array}{llll}
SGp(G(K)\bar{\wedge}L, {\bf H})&\cong& S_{\ast}(L, SGp(G(K), {\bf H}))\\
& \cong& S_{\ast}(L, S_{\ast}(K, \bar{W}{\bf H})) \\
& \cong& S(K\wedge L, \bar{W}{\bf H})\\
& \cong & SGp(G(K\wedge L), {\bf H})
\end{array}
$$
thus  $G(K)\bar{\wedge}L\cong G(K\wedge L).$ As Curtis notes (\cite{curtis} p. 137) 
$\bar{W}K(\pi, 0)$ is a minimal complex for $K(\pi,1)$ so taking $K = \bar{W}K(\pi,0)$ 
we get $\Omega\sum K(\pi, 1)$ has as model $G(\bar{W}(K(\pi, 0))\wedge S^1)$ and hence 
$G(\bar{W}(K(\pi, 0)))\bar{\wedge}S^1.$ By Lemma, 4.4 given the  weak homotopy equivalence 
$G\bar{W}(K(\pi, 0))\longrightarrow K(\pi,0),$ the result follows. \end{pf}

This implies that, like Jie Wu \cite{wu}, we can take a simple model for $\Omega\sum K(\pi, 1).$ 
First we introduce notation for $S^1.$ We write $S_0^1 =\ast, \ S_1^1 =\{\sigma,\ast\}, \ 
S_2^1 =\{ x_0, x_1, \ast \}$ where $x_0 = s_1\sigma, \ x_1 = s_0\sigma$ and in general 
$S_{n+1}^1 =\{x_0,\ldots,x_n,\ast\}$ where 
$x_i =s_n\ldots s_{i+1}s_{i-1}\dots s_0\sigma,$ $0\leq i \leq n.$

For simplicity we  write $G= K(\pi,0)$ and make no distinction between simplicies in 
different dimensions unless confusion might arise. This then gives 
$$
\begin{array}{llll}
(G\bar\wedge S^1)_0 &=& 1, \qquad\text{the trivial group}\\
(G\bar\wedge S^1)_1 &\cong& \pi, \\
(G\bar\wedge S^1)_1 &\cong&\pi\ast\pi, \qquad\text{the free product of two copies of $\pi$}\\
(G\bar\wedge S^1)_3 &\cong&\pi\ast\pi\ast\pi \qquad\text{and so on}.
\end{array}
$$
We write $g\bar{\wedge}x$ for the $x$-indexed copy of $g\in \pi$ in the coproduct 
$\coprod\{(\pi)_x : x\in S_n^1\setminus\{{\ast}\}\}$ so the only relations we have are of the form
$$
(g{g'}\bar{\wedge} x) = (g\bar{\wedge} x)({g'}\bar{\wedge} x).
$$

We next analyse $N(G\bar{\wedge}S^1)$ in low dimensions. For simplicity we will write 
${\bf H}$ instead of $G\bar{\wedge}S^1.$ Of course $NH_0 =1,\ NH_1 = \pi.$

By Theorem A, $NH_2$ is generated by all $F_{(0)(1)}(g\bar{\wedge}\sigma, h\bar{\wedge}\sigma),$ 
$g,\ h \in \pi:$
$$
F_{(0)(1)}(g\bar{\wedge}\sigma, h\bar{\wedge}\sigma)= 
(g\bar{\wedge}x_1)(h\bar{\wedge}x_0)(g^{-1}\bar{\wedge}x_1)(gh^{-1}g^{-1}\bar{\wedge}x_0).
$$

In fact although Theorem A gives these as normal generators, it is clear that these 
are generators since conjugates of them are expressible as product of other
terms of the same form. For instance
$$
{}^{k\bar{\wedge}x_1}F_{(0)(1)}(g\bar{\wedge}\sigma, h\bar{\wedge}\sigma)= 
F_{(0)(1)}(kg\bar{\wedge}\sigma, h\bar{\wedge}\sigma)
F_{(0)(1)}(kg\bar{\wedge}\sigma, ghg^{-1}\bar{\wedge}\sigma)
$$
and a similar expression can be found for conjugation by $k\bar{\wedge}x_0.$

Brown and Loday \cite{bl1} calculated $\pi_3(\sum K(\pi, 1))$ using a van Kampen
theorem for cat$^2$-groups. This led to an expression for this group as being isomorphic to
$$
J_2(\pi) = \mbox{Ker}(\kappa: \pi\otimes\pi\longrightarrow\pi).
$$
We refer to the paper \cite{bjr} Brown, Johnson and Robertson for some details on the non-abelian tensor
product of groups and to Ellis \cite{ellis2} for more on the representation of homotopy types by
crossed squares and cat$^2$-groups. Here it will suffice to say that if $G,$ $H$ are groups that act 
on themselves by conjugation and on each other in such a way that
$$
{}^{{}^{g}h}{g'}={}^{ghg^{-1}}{g'}\hspace{2cm}{}^{{}^{h}g}{h'}={}^{hgh^{-1}}{h'}
$$
(see \cite{bjr}), the tensor product $G\otimes H$ is the group generated by symbols $g\otimes h$
with relations
$$
\begin{array}{lcr}
g{g'}\otimes h &=& ({}^{g}{g'}\otimes {}^{g}h)(g\otimes h)\cr
g\otimes h{h'} &=& (g\otimes h)\otimes({}^{h}g\otimes {}^{h}{h'})
\end{array}
$$
for all $g,{g'}\in G,$ $h,{h'}\in H.$ We will need this only when $G = H.$, then there is a map 
$$
\kappa: G\otimes G\longrightarrow G
$$
given by $\kappa(g\otimes h)=[g,~h].$ This gives a homomorphism since the above relations are abstract
versions of the usual commutator identities.

Using a combination of Ellis's work in \cite{ellis2} and the second author's description of the
crossed $n$-cube associated to a simplicial group in \cite{porter}, it is clear that the expression
for the tensor product should be closely linked to one for $NH_2/d_3NH_3$ and with that in mind we set
for $g,~h\in\pi$
$$
g\bar{\otimes}h = [g\bar{\wedge}x_0,~(h\bar{\wedge}x_0)(h^{-1}\bar{\wedge}x_1)]~d_3NH_3.
$$
This `mysterious' formula will be more fully explained in another paper where the relationship 
with crossed squares and Ellis's work will be given in detail. For  the moment the reader is asked 
merely to accept the left hand side as a shorthand for the right hand side.
\begin{prop}
The symbols $g\bar{\otimes} h$ generate $NH_2/d_3NH_3$ and satisfy the tensor product relations above.
The boundary homomorphism
$$
{d_2'}:\frac{NH_2}{d_3NH_3}\longrightarrow NH_1
$$
sends $g\bar{\otimes} h$ to $[g,~h]\bar{\wedge}\sigma$ in $NH_1.$
\end{prop}
\begin{pf}
Direct calculation shows
$$
g\bar{\otimes} h = F_{(0)(1)}(h\bar{\wedge}\sigma, g\bar{\wedge}\sigma)~d_3NH_3
$$
which clearly generate $NH_2/d_3NH_3$ by Theorem A and the commutator above. 
The tensor product relations are again directly verifiable (eg. using the description of the 
$h$-map given in \cite{porter}), and it is immediate that
$$
d_2[g\bar{\wedge}x_0,~(h\bar{\wedge}x_0)(h^{-1}\bar{\wedge}x_1)]= [g,~h].
$$\hfill\end{pf}

The next term in the Moore complex $NH_3$ has 6 different types of generator and so there are 6
$\underline{\mbox{known}}$ types of relation in $NH_2/d_3NH_3.$ Presumably it is possible to give 
a direct proof that 
$$
\frac{NH_2}{d_3NH_3}\cong\pi\otimes\pi
$$
using these relation, but we have so far not managed to find one. A closer analysis of the relationship
of the Peiffer pairings with the structure of crossed squares gives this isomorphism by a universal 
argument, but requires other techniques and so will be postponed to a later paper. Knowledge of the 6
types of generator for $NH_3$ hopefully will allows a detailed calculation of $NH_3/d_4NH_4$ in a similar
manner.


{\noindent  A. MUTLU \hspace{7.3cm} T. PORTER} \\
{\it Department of Mathematics \hspace{4.4cm}  School of Mathematics \\                        
Faculty of Science  \hspace{6.cm}  Dean Street\\                               
University of Celal Bayar  \hspace{4.8cm} University of Wales, Bangor\\                          
Manisa, TURKEY  \hspace{6.cm} Gwynedd, LL57 1UT, UK}\\ 
e-Mail: amutlu@spil.bayar.edu.tr \hspace{3.4cm}e-Mail: t.porter@bangor.ac.uk\\

\end{document}